\documentclass[12pt]{article}
\usepackage{amsmath,amssymb,amsthm, enumerate, fullpage,graphicx, hyperref,verbatim,ytableau}

\newtheorem{theorem}{Theorem}
\newtheorem{lemma}[theorem]{Lemma}

\newtheorem{proposition}[theorem]{Proposition}

\AtBeginDocument{}

\author{ \Large{ Matt Davis } \\ \\ Department of Mathematics and Computer Science\\ Muskingum University,  10 College Dr., New Concord, OH, 43762 \\ mattd@muskingum.edu }
\title{A bijection between all Shi regions and core partitions}
\date{}

\begin{document}

\maketitle

\section{Introduction}

The $m$-Shi arrangement is a hyperplane arrangement related to a number of areas of study in combinatorics and representation theory, including parking functions and $n$-core partitions. There are multiple bijections between the set of all regions of the arrangement and parking functions (\cite{AL}, \cite{PS}), and the Athanasiadis-Linusson bijection in particular is interesting in part because the natural symmetric group action on parking functions reflects the geometry of the Shi arrangment.

In \cite{FV1} (resp. \cite{FV2}), Fishel and Vazirani showed a bijection between the dominant (resp. bounded dominant) chambers of the $m$-Shi arrangement and partitions which are simultaneously $n$-core and $(mn+1)$-core (resp. $(mn-1)$-core). However, this bijection does not extend outside the dominant chamber. The goal of this paper is to extend the results of Fishel-Vazirani by giving a bijection between the set of all chambers of the $m$-Shi arrangement and a certain subset of $n$-core partitions. Moreover, like the original, it adapts in a natural way to a bijection on the bounded chambers of the arrangement. In addition, it has the same $S_{n}$-structure as parking functions in the Athanasiadis-Linusson bijection. Thus we close the loop from the Fishel-Vazirani result to the more classical indexings of Shi regions, and bring the $S_{n}$ action to the partition picture. As part of the proof, we will establish an independently interesting characterization of the minimal chambers that occur in the $S_{n}$-orbit of a minimal chamber of the $m$-Shi arrangement.

\section{The Affine Symmetric Group and Alcoves}

We being with a brief survey of known results on the relationship between the $m$-Shi arrangement, $n$-core partitions, and the affine symmetric group. Our notation in this section closely models that of \cite{FV1}, which contains several results which will be important in what follows. 

\subsection{The affine symmetric group} The affine symmetric group $\widetilde{S}_{n}$ is defined to be the set of all permutations $\phi$ of $\mathbb{Z}$ satisfying the condition that, for any integer $m$, \[ \phi(m+n) = \phi(m)+n\textrm{, and } \sum_{i=0}^{n-1} \phi(i) = \sum_{i=0}^{n-1} i.\] To define an element $\phi \in \widetilde{S}_{n}$, it is sufficient to specify $\phi(0), \phi(1), \cdots \phi(n-1)$, a set of numbers that form a transversal of the congruence classes of $\mathbb{Z}$ mod $n$, and sum to $\binom{n}{2}$. We will refer to such an $n$-tuple $(\phi(0), \phi(1), \ldots \phi(n-1))$ as the $n$-\textit{window} of $\phi$.

This group is generated by the elements $s_{0}, s_{1}, \cdots s_{n-1}$, where $s_{0}$ has the $n$-window $(-1, 1, 2, \ldots, n-2,n)$, and for $1 \leq i \leq n-1$, $s_{i}$ has $n$-window $(0,1,\ldots i-2, i, i-1,i+1, \ldots n-1)$. With these generators, $\widetilde{S}_{n}$ for $n>2$ has presentation: \[ s_{i}^{2} = e \] \[ s_{i}s_{j} = s_{j}s_{i} \textrm{ if } j \neq i \pm 1 \textrm{ mod } n\] \[ s_{i}s_{j}s_{i} = s_{j}s_{i}s_{j} \textrm{ if } j = i \pm 1 \textrm{ mod } n. \] We will always assume that $n>2$ in what follows. 

We also note that $G = \langle s_{1}, \ldots s_{n-1} \rangle$ is a subgroup of $\widetilde{S}_{n}$ isomorphic to the finite symmetric group $S_{n}$. There is a homomorphism $f_{n}: \widetilde{S}_{n} \rightarrow G$ given by sending $s_{0}$ to $s_{\theta} = s_{1}s_{2}\ldots s_{n-2}s_{n-1}s_{n-2}\ldots s_{2}s_{1}$ and fixing each other $s_{i}$. We will use $G$ and $f_{n}$ in this paper only to mean this specific group and homomorphism.

\subsection{Root data} We will need to understand the root data associated to the affine symmetric group $\widetilde{S}_{n}$. We let $\epsilon_{1}, \ldots \epsilon_{n}$ be the standard basis for $\mathbb{R}^{n}$, and let $\langle \cdot | \cdot \rangle$ be the usual inner product. Let $V = \{ a_{1} \epsilon_{1} + \ldots + a_{n} \epsilon_{n} \, | \, a_{1} + \ldots + a_{n} = 0\}$, and let $\alpha_{i} = \epsilon_{i} - \epsilon_{i+1}$. Then $\Pi = \{ \alpha_{i} \, | \, 1 \leq i \leq n-1\}$ is a basis for $V$, and $\Pi$ is the set of simple roots of the type $A_{n-1}$ root system. We identify $Q = \bigoplus_{i}^{n-1} \mathbb{Z}\alpha_{i}$ with the root lattice of this system. The set of all roots is $\Delta = \{\epsilon_{i} - \epsilon_{j} , \, | \, 1 \leq i,j, \leq n, \, i \neq j \}$, and $\theta = \alpha_{1} + \ldots + \alpha_{n} = \epsilon_{1} - \epsilon_{n}$ is the highest root.  The set $\Delta^{+}$ of positive roots consists of those roots $\epsilon_{i} - \epsilon_{j}$ where $i < j$. We let $\Delta^{-} = - \Delta^{+}$ and write $\alpha >0$ (resp. $\alpha < 0$) if $\alpha \in \Delta^{+}$ (resp. $\alpha \in \Delta^{-}$). As usual, for any root $\alpha \in \Delta^{+}$ we will use $s_{\alpha}$ to denote the reflection associated to $\alpha$ in the Weyl group $S_{n}$, with $s_{i} = s_{\alpha_{i}}$ and $s_{\theta} = s_{1}s_{2} \ldots s_{n-1} \ldots s_{1}$ as above.

We must also consider the affine root system $\widetilde{A}_{n-1}$, built by including the ``imaginary'' root $\delta$ in the root system of type $A_{n-1}$. We let $\widetilde{\Delta} = \{k \delta + \alpha \, | \, \alpha \in \Delta, k \in \mathbb{Z}\}$ with $\widetilde{\Delta}^{+} = \{k\delta + \alpha \, | \, k \in \mathbb{Z}_{\geq 0}, \alpha \in \Delta^{+} \textrm{ OR } k \in \mathbb{Z} > 0, \alpha \in \Delta^{-}\}$, and $\widetilde{\Delta}^{-} = - \widetilde{\Delta}^{+}$. We let $\alpha_{0} = \delta - \theta$, and the simple roots in $\widetilde{A}_{n-1}$ are $\widetilde{\Pi} = \Pi \cup \{ \alpha_{0}\}$. We again write $\alpha > 0$ (resp. $\alpha<  0$) to mean $\alpha \in \widetilde{\Delta}^{+}$ (resp. $\alpha \in \widetilde{\Delta}^{-}$.) 

For any $\alpha \in \Delta$ and $k \in \mathbb{Z}$, we define the hyperplane \[ H_{\alpha,k} = \{ v \in V \, | \, \langle v | \alpha \rangle = k\}.\] Note that $H_{-\alpha,-k} = H_{\alpha,k}$. For any hyperplane $H_{\alpha,k}$, we consider the (closed) half spaces $H_{\alpha,k}^{+} = \{ v \in V \, | \, \langle v | \alpha \rangle \geq k\}$ and $H_{\alpha,k}^{-} = \{ v \in V \, | \, \langle v | \alpha \rangle \leq k\}$.

\subsection{Hyperplane Arrangements} We will be mainly interested in two particular hyperplane arrangements: the \textit{Coxeter arrangement} $\mathcal{A}_{n}$, consisting of all the hyperplanes $\{H_{\alpha,k} \, | \, \alpha \in \Delta^{+}, k \in \mathbb{Z}\}$, and the \textit{$m$-Shi arrangement} $\mathcal{S}_{n,m}$, consisting of the hyperplanes $\{ H_{\alpha, k} \, | \, \alpha \in \Delta^{+} \textrm{ and }  -m < k \leq m \}$. For a hyperplane arrangement $\mathcal{H}$, we are interested in the connected components of $V \setminus \bigcup_{H \in \mathcal{H}} H$, the complement of the hyperplanes in $\mathcal{H}$. Each such connected component is called a \textit{region} of the hyperplane arrangement, and can be written as the intersection of a collection of half-planes of hyperplanes in $\mathcal{H}$. A region is said to be \textit{dominant} if it is a subset of $\bigcap_{i=1}^{n} H_{\alpha_{i},0}^{+}$. We call $\bigcap_{i=1}^{n} H_{\alpha_{i},0}^{+}$ the \textit{dominant chamber}, and the elements of its $G$-orbit are called \textit{chambers} or \textit{Weyl chambers}.

The regions of $\mathcal{A}_{n}$ are referred to as \textit{alcoves}, and $A_{0} =  \left(\bigcap_{\alpha \in \Pi^{+} } H_{\alpha,0}^{+}\right) \cap H_{\theta,1}^{-}$ is the \textit{fundamental alcove}. Note that since $\mathcal{S}_{n,m} \subseteq \mathcal{A}_{n}$, the regions of $\mathcal{S}_{n,m}$ are (up to a measure 0 set) unions of alcoves. As such, the geometry of $\mathcal{A}_{n}$ is fundamental to understanding $\mathcal{S}_{n,m}$. In turn, the affine symmetric group is intimately tied to the set of alcoves of $\mathcal{A}_{n}$; we summarize the relationship here.

We let $\widetilde{S}_{n}$ act on $V$ via: \[ s_{i} \cdot (a_{1}, \ldots, a_{n}) = (a_{1}, \ldots, a_{i-1}, a_{i+1}, a_{i}, a_{i+2}, \ldots, a_{n}) \textrm{ for } i \neq 0 \]
\[ s_{0} \cdot (a_{1}, \ldots, a_{n}) = (a_{n} + 1, a_{2}, \ldots a_{n-1}, a_{1} - 1).\] Then for $i \neq 0$, $s_{i}$ acts as reflection through $H_{\alpha_{i},0}$, and $s_{0}$ reflects through $H_{\theta,1}$. This action preserves the lattice $Q$ and acts freely and transitively on the set of alcoves. Thus we can identify each element $w \in \widetilde{S}_{n}$ with the alcove $w \cdot A_{0}$, and this correspondence is a bijection.

Note that $G$ preserves $\langle | \rangle$ on $V$ (although $\widetilde{S}_{n}$ does not). Let $w \in \widetilde{S}_{n}$. Then if $y$ is a minimal-length representative of the coset $Gw$, then either $y=e$ or any reduced word for $y$ must have $s_{0}$ as its leftmost letter.  We note here (see \cite{L}, for example) that if $w \in \widetilde{S}_{n}$ and $y$ is a minimal-length coset representative of $Gw$, then the window for $y^{-1}$ is sorted in ascending order, and $yA_{0}$ is in the dominant chamber. We will abuse notation slightly and use $\widetilde{S}_{n}/G$ to denote the set of minimal-length right coset representatives of $G$ in $\widetilde{S}_{n}$.

The combinatorics of the alcoves of $\mathcal{A}_{n}$ and $\widetilde{S}_{n}$ are closely related, but here we will only point out the fact that the length of an element $w \in \widetilde{S}_{n}$ is exactly the number of distinct hyperplanes in $\mathcal{A}_{n}$ separating $wA_{0}$ from the origin.

The boundary of any alcove of $A_{n}$ is a union of subsets of various hyperplanes $H_{\alpha,k}$, which we call \textit{facets}. Each facet of an alcove is in the $\widetilde{S}_{n}$-orbit of exactly one facet of $A_{0}$ (\cite{S1}). If $H_{\alpha,k}$ is the hyperplane containing a facet of an alcove, we call it a \textit{wall} of the alcove. Since the walls of $A_{0}$ are $H_{\alpha_{i},0}$ for $1 \leq i \leq n-1$ or $H_{\theta,1}$, we can label each of these facets with the corresponding $i$ or 0 for $H_{\theta,1}$. Applying the same label to each facet in the orbit of the facets of $A_{0}$ applies a label to each facet of each alcove of $A_{n}$. 

Shi proved (\cite{S1}) that each region of $\widetilde{S}_{n,m}$ contains a unique alcove $wA_{0}$ such that $\ell(w)$ is minimal. We call such an alcove \textit{$m$-minimal}, and note that studying the regions of $\mathcal{S}_{n,m}$ is in many ways equivalent to studying the $m$-minimal alcoves. Moreover, any bounded region of $\mathcal{S}_{n,m}$ has a unique alcove $wA_{0}$ such that $\ell(w)$ is maximal, which we call an \textit{$m$-maximal region.}

Fishel-Vazirani gave the following characterization of the walls of an alcove.

\begin{lemma} \label{lem:mainFV} (\cite{FV1}, Proposition 4.1) For $w \in \widetilde{S}_{m}$ and $0 \leq i \leq n-1$, $wA_{0} \subseteq H_{\alpha,k}^{+}$ and $ws_{i}A_{0} \subseteq H_{\alpha,k}^{-}$ if and only if $w(\alpha_{i}) = \alpha-k\delta$. \end{lemma}

The first condition here says that the alcove $wA_{0}$ has $H_{\alpha,k}$ as one of its walls, with 
label $s_{i}$, and the alcove $ws_{i}A_{0}$ on the other side of the wall is on the negative side of $H_{\alpha,k}$. Note that \cite{FV1} only proves the forward implication, but the reverse is implicit. If $w(\alpha_{i}) = \alpha-k\delta$, then we note that $w\mathcal{A}_{0}$ always has some facet with label $i$, so $ws_{i}\mathcal{A}_{0}$ shares a facet with it. Then if $wA_{0} \subseteq H_{\alpha',k'}^{+}$ and $ws_{i}A_{0} \subseteq H_{\alpha',k'}^{-}$ (possibly for $\alpha' < 0$), we must have $w(\alpha_{i}) = \alpha' - k'\delta$, so $k' = k$ and $\alpha' = \alpha$.

 A \textit{floor} of an alcove is a wall of that alcove that does not go through the origin, but separates the alcove from the origin. A \textit{ceiling} of an alcove is a wall of that alcove so that the alcove and origin are on the same side of the wall.

\section{n-cores and abaci}

We now summarize some standard facts about partitions - see \cite{JK} for example, for details. A \textit{partition} $\lambda$ of a positive integer $n$ is an infinite sequence $(\lambda_{1}, \lambda_{2}, \cdots)$ of non-negative integers, written in non-increasing order, that sum to $n$. For example, $(5,4,2,1, 0,0, \cdots)$ and $(10,1,1,0,0,\cdots)$ are partitions of 12. We will usually write a partition as a finite list consisting only of the non-zero entries of $\lambda$. We will also typically identify a partition with its Young diagram - a configuration of boxes arranged in a grid, where the first row has $\lambda_{1}$ boxes, the second has $\lambda_{2}$ boxes, etc.

\begin{figure}[ht]
\centering
\ytableausetup{centertableaux}
\begin{ytableau}
0 & 1 & 2 & 3 & 4 \\
-1 & 0 & 1  \\
-2  \\
-3
\end{ytableau}
\caption{The partition $(5,3,1,1)$, with contents}
\end{figure}

The boxes in $\lambda$ will be referred to by their coordinates, so that box $(i,j)$ is in the $j$th row and $i$th column, counting from the top and left sides of the diagram. We say that a box $b$ in a Young diagram $\lambda$ is \textit{removable} for $\lambda$ if it is in $\lambda$, but removing it from $\lambda$ still yields a valid Young diagram. We say that a box $b$ is \textit{addable} for $\lambda$ if it is not in $\lambda$, but adding it to $\lambda$ would still yield a valid Young diagram. Each box in the grid is assigned a \textit{content}. Specifically, the box in row $j$ and column $i$ has content $i-j$.

Each box in a Young diagram also has associated to it a \textit{hook length}, defined as the the number of boxes directly to the right of the box in its row plus the number of boxes directly below the box in its column, plus 1.

\begin{figure}[ht]
\centering
\ytableausetup{centertableaux}
\begin{ytableau}
8 & 5 & 4 & 2 & 1 \\
5 & 2 & 1 \\
2 \\
1
\end{ytableau}
\caption{The hook lengths for $(5,3,1,1)$}
\end{figure}

If a given box $b$ has hook length $s$, then there will be a connected chain of boxes along the lower-right edge of the Young diagram, starting at the end of the row containing $b$, and ending at the bottom of the column containing $b$, that has exactly $s$ boxes.

\begin{figure}[ht]
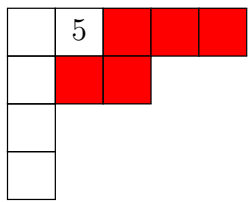

\centering
\ytableausetup{centertableaux}
\begin{ytableau}
\, & 5 & *(red) & *(red) & *(red) \\
\, & *(red) & *(red) \\
\,  \\
\,
\end{ytableau}
\caption{The rim hook for box $(2,1)$}
\end{figure}

This chain of boxes is the \textit{rim hook} associated to $b$. Removing the rim hook associated to a given box will yield a valid Young diagram. If we repeatedly remove rim hooks of length $n$ from a Young diagram $\lambda$ until none remain, the result will be independent of the order in which the rim hooks are removed. In this case, we call the resulting partition $\lambda'$ the $n$-core of $\lambda$. If a Young diagram $\lambda$ has no hook lengths equal to $n$, we say that $\lambda$ is $n$-core, since the $n$-core of $\lambda$ is $\lambda$ itself. A partition $\lambda$ is $n$-core if and only if none of its hook lengths are \textit{divisible} by $n$. We say that a partition $\lambda$ is $(n,t)$-core if it is simultaneously $n$-core and $t$-core.

In \cite{L}, Lascoux described a way for $\widetilde{S}_{n}$ to act on Young diagrams of $n$-core partitions by adding or removing boxes of a given content. Specifically, if $\lambda$ is a partition, then it cannot have both addable and removable boxes with content congruent to $i$ mod $n$. (Due to the importance of this action, in our treatment of $n$-cores, we will always reduce the contents of boxes mod $n$.) Then $s_{i}$ acts on $\lambda$ by

\[ s_{i} \cdot_{1} \lambda = \begin{cases} \lambda \cup \{ \textrm{All addable boxes of content $i$ mod $n$} \} & \textrm{ If $\lambda$ has such addable boxes} \\ \lambda \setminus \{ \textrm{All removable boxes of content $i$ mod $n$} \} & \textrm{ If $\lambda$ has such removable boxes} \\ \lambda & \textrm{ Otherwise } \end{cases} .\]

This is the description of this action originally given in \cite{L}, although \cite{FV1} has more of the details worked out. We will refer to this as the level 1 action, to contrast the level $t$ action below.

There are a number of combinatorial indexings of $n$-core partitions that we will move between often in this paper. We will begin by introducing all of these indexings and a description of the necessary actions of $\widetilde{S}_{n}$ on each one.

\textbf{Balanced $n$-abaci} The first combinatorial indexing of partitions which is particularly useful for $n$-cores is the notion of an $n$-abacus (\cite{JK} is a standard reference). Imagine the integers arranged in an infinite grid containing $n$ columns, marked with the numbers $0,1,2, \cdots n-1$. We label the rows by integers, so that row 0 contains $0,1,2, \cdots n-1$, row 1 contains $n,n+1, \ldots 2n-1$,  row -1 contains $-n,-n+1, \ldots -1$, etc. An abacus diagram is merely a set of integers whose complement is bounded below, which we draw on the infinite grid by circling each number in the set. The columns of this diagram are called the \textit{runners} of the abacus, and the circled numbers are called \textit{beads}.

From a partition $\lambda$, we can obtain a particular abacus diagram which we will call the \textit{positive $n$-abacus} of $\lambda$. The positive $n$-abacus of a partition $\lambda$ is the diagram obtained from this grid by circling each number which is the hook length of a box in the first column of the Young diagram of $\lambda$, as well as all the negative integers. Then a partition is $n$-core if and only its positive $n$-abacus has no positive beads on the 0 runner, and for each $i$ from $1$ to $n-1$, if runner $i$ has any positive beads, they are in positions $i, i+n, i+2n, \ldots i+kn$ for some nonnegative integer $k$. In other words, the beads on the abacus must be pushed as far up the runners as possible, with no positive beads on the 0 runner.

We say that two abacus diagrams are equivalent to each other if one can be obtained from the other by adding a constant to each bead in the diagram. Any abacus diagram that is equivalent to the positive $n$-abacus of $\lambda$ will be referred to as an $n$-abacus for $\lambda$.

\begin{figure}[ht]
    \centering
    \includegraphics[width=0.5\linewidth]{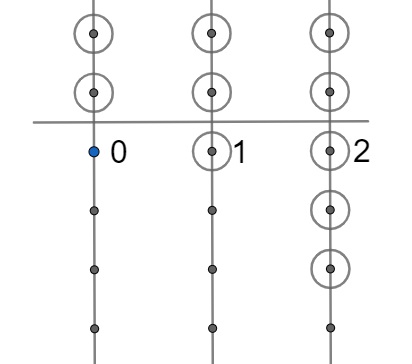}
    \caption{The positive 3-abacus diagram for $(5,3,1,1)$}
    \label{fig:Abacus}
\end{figure}

We define the \textit{balance number} of an abacus diagram for an $n$-core $\lambda$ to be the sum over the runners of the diagram of the row number of the first non-circled number on that runner, and say that a diagram is balanced if its balance number is 0. 

\begin{figure}[ht]
    \centering
    \includegraphics[width=0.5\linewidth]{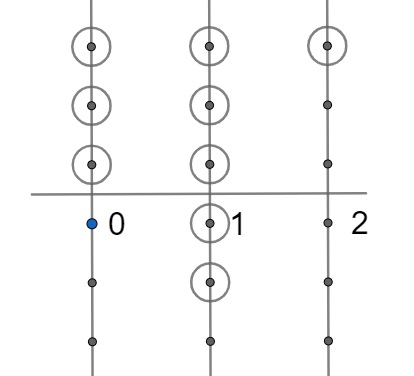}
    \caption{The balanced 3-abacus diagram for $(5,3,1,1)$}
    \label{fig:Abacus2}
\end{figure}

Then every $n$-core partition has an abacus diagram equivalent to   a unique balanced abacus diagram, which we call \textit{the} balanced $n$-abacus of $\lambda$. This correspondence gives a bijection between balanced $n$-abacus diagrams meeting the ``flushness'' property described above and $n$-core partitions. (Note the positive part of the 0 runner may no longer be empty in the balanced abacus.) The balanced $n$-abacus is particularly useful because the action of $\widetilde{S}_{n}$ on the set of balanced $n$-abaci is easily described. Given a balanced $n$-abacus $X$ for $\lambda$, we define $w_{i} \cdot X$ to be $X$ with runners $i-1$ and $i$ swapped if $1 \leq i \leq n-1$, and $w_{0} \cdot X$ to be the diagram obtained from $X$ by swapping runners $0$ and $n-1$, and then adding one bead to runner 0 and removing the last bead from runner $n-1$. The result will be the balanced abacus of an $n$-core partition, and \cite{BJV} shows that this action is equivalent to the action on $n$-core partitions described above.

\textbf{$n$-vectors} Because of the ``flushness'' condition, the balanced $n$-abacus of an $n$-core partition is completely determined by the row number of the first non-circled bead  on each runner. We can encode this as a vector with $n$ components, the $i$th entry of which is the row number of the first non-circled number on that runner. This is the `$n$-vector' construction of Garvan, Kim, and Stanton (\cite{GKS}). Since the abacus diagram is balanced, the entries of the resulting vector will sum to 0. Thus this correspondence gives a bijection between $n$-core partitions and $n$-dimensional vectors with integer entries that sum to zero. In \cite{BJV} it is shown that this is an $\widetilde{S}_{n}$-equivariant bijection between $n$-core partitions and $Q$.

Moreover, \cite{FV1} shows that this map induces a bijection between the set of partitions which are $(n,mn+1)$-core and the dominant $m$-minimal alcoves. Specifically, if $\lambda$ is an $(n,mn+1)$-core partition equal to $w^{-1} \cdot \emptyset$, where $w \in \widetilde{S}_{n}/G$, then $wA_{0}$ is an $m$-minimal alcove. 

\textbf{$n$-sets and $n$-windows} Consider the $n$-vector $v = (a_{1}, \ldots, a_{n})$ of an $n$-core partition $\lambda$, and let $S(\lambda) = \{na_{1}, na_{2} + 1, \ldots, na_{n} + n-1\}$. Since $v$ encodes the rows of the first non-circled beads on each runner of the abacus, the (unique) element of $S(\lambda)$ that is congruent to $i$ mod $n$ is the label of the first non-circled bead on runner $i$. Also, since $\sum a_{i} = 0$, we notice that the sum of the elements of $S(\lambda)$ is $\binom{n}{2}$. The tuple $S(\lambda)$ is called the $n$-set of $\lambda$, and the resulting correspondence is a bijection between transversals of $\mathbb{Z}/n\mathbb{Z}$ with sum equal to $\binom{n}{2}$ and $n$-core partitions. Since $\widetilde{S}_{n}$ acts on $\mathbb{Z}$, it also acts on $n$-sets by acting on each element of the set. It is routine to check that the action on the $n$-set of a partition corresponds to the action on its abacus.

We will use $S(\lambda)_{i}$ to denote the element of $S(\lambda)$ which is congruent to $i$ mod $n$.  Finally, we let $X(\lambda)$ be the $n$-tuple consisting of the elements of $S(\lambda)$ sorted in increasing order, which we call the $n$-window for $\lambda$. Similarly to above, we will use $X(\lambda)_{i}$ to denote the $i$th entry of $X(\lambda)$, i.e. the $i$th smallest entry in $S(\lambda)$. Then if $\lambda = w^{-1}\emptyset$ where $w \in \widetilde{S}_{n}/G$, we have that $X(\lambda)$ is the $n$-window of $w^{-1}$. We note that $n$-windows give yet another parametrization for $n$-cores.

\section{The level t action and n-cores}

We turn now to Fayers' ``level $t$'' action of $\widetilde{S}_{n}$ on $\mathbb{Z}$. Specifically, in \cite{F1} and \cite{F2}, for any $t$ relatively prime to $n$ (which we assume throughout), Fayers defined an action which we denote as $\ast_{t}$:

\[ s_{i} \ast_{t} j = \begin{cases} j+t & \textrm{ if } j \equiv (i-1)t - n \circ t  \textrm{ mod } n \\  j-t & \textrm{ if } j \equiv it - n \circ t \textrm{ mod } n \\ j & \textrm{ otherwise } \end{cases} ,\] where $n \circ t = \frac{1}{2}(n-1)(t-1)$. We will not use this precise action, but we restate it here to point out why Fayers' results still apply in our case. Instead, we focus on two variants of this action which are specific to the case $t \equiv \pm 1$ mod $n$, denoted $\cdot_{t}$. First, if $t = mn+1$, we define:

\[ s_{i} \cdot_{t} j = \begin{cases} j+t & \textrm{ if } j \equiv (i-1) \textrm{ mod } n \\  j-t & \textrm{ if } j \equiv i \textrm{ mod } n \\j & \textrm{ otherwise } \end{cases} .\]

Secondly, if $t = mn-1$, then we define:

\[ s_{i} \cdot_{t} j = \begin{cases} j-t & \textrm{ if } j \equiv (i-1) \textrm{ mod } n \\  j+t & \textrm{ if } j \equiv i \textrm{ mod } n \\j & \textrm{ otherwise } \end{cases} .\]

We use these particular definitions because in either case, the level $t$ action of $s_{i}$ will swap the $i-1$ and $i$ congruence classes setwise, by shifting by $t$ in the appropriate direction. A similar definition could be given done for any value of $t$ relatively prime to $n$, but the choice of which congruence classes to swap may start to matter, and it is unclear if such an action has a direct application. However, in our case, tailoring the level $mn+1$ and $mn-1$ actions in this way will let us give a consistent description of both of our bijections. (Note that we also have the nice property that the level 1 action is the usual action of $\widetilde{S}_{n}$, which we will denote as $\cdot_{1}$ from now on. )

Now, for $t=mn+1$, the actions $\ast_{t}$ and $\cdot_{t}$ differ only by an index shift by $n \circ t$, an automorphism of $\widetilde{S}_{n}$. Similarly, for $t = mn-1$, the actions differ by shifting the index by $n \circ t$ and negating (mod $n$). Thus, the orbit of any integer or set of integers will be the same under either $\ast_{t}$ or $\cdot_{t}$, meaning several important results from \cite{F1} still hold. (We will refer to simply ``the orbit of the level $t$ action", since it is the same regardless of whether we look at $\cdot_{t}$ or $\ast_{t}$. For results specific to $t = mn \pm 1$, we will specify the value of $t$.) The proofs of the other results we need from \cite{F1} are virtually unchanged, and we  summarize all those results here.

The level $t$ action is important for us in the way that it acts on $n$-windows, $n$-cores, and their $n$-sets. In \cite{F1}, Fayers described $\ast_{t}$ acting on $n$-cores, which operates in terms of adding or removing rim hooks of legnth $t$ with a particular pattern of contents - a direct generalization of Lascoux's $\widetilde{S}_{n}$ action. On the level of $n$-sets, the level $t$ action on $\mathbb{Z}$ can be applied directly to sets of integers as well. Fayers (\cite{F1}) shows that if $\lambda$ is an $n$-core, then for any $w \in \widetilde{S}_{n}$, $w \ast_{t} S(\lambda)$ is the $n$-set of another $n$-core partition. Now, as in \cite{F2} we let $C_{n}$ be the set of all $n$-cores, and define $C_{n}^{(N)}$ to be the set of all $n$-cores $\lambda$ satisfying $|x -y| < nN$ for all $x,y \in S(\lambda)$. Fayers (\cite{F2}, Lemma 3.16) shows that each $(s,t)$ core is an element of $C_{n}^{(t)}$. We also give an equivalence relation $\equiv_{N}$ on $n$-cores where $\lambda_{1} \equiv_{N} \lambda_{2}$ if there is a bijection $\phi: S(\lambda_{1}) \rightarrow S(\lambda_{2})$ with $\phi(x) \equiv x$ mod $nN$ for all $x \in S(\lambda_{1})$.

We should note here that $C_{n}^{(mn+1)}$ has essentially already appeared as an indexing set for chambers of the Shi arrangement. The group $\mathbb{Z}^{n}_{n+1}$ of $n$-tuples of integers mod $n+1$ has a subgroup $H$ of order $n+1$ generated by $(1,1,\ldots,1)$, so the quotient $\mathbb{Z}^{n}_{n+1}/H$ has $(n+1)^{n-1}$ elements. (This set appears in \cite{A1}, \cite{AL}, and \cite{S}, among many others.) But each coset has a unique representative with a 0 in the first entry, and $(0,a_{1},\ldots a_{n-1})$ can be associated to the positive abacus with $a_{i}$ (flush) beads on the runner labeled $i$. If all the $a_{i}$ are less than $n+1$, then that abacus represents a partition in $C_{n}^{(n+1)}$. These cosets also each contain a single representative that describes a parking function, and the Athanasiadis-Linusson bijection (though not the Pak-Stanley) from regions to parking functions is $S_{n}$-equivariant with the normal permutation action on $\mathbb{Z}^{n}_{n+1}/H$. Interestingly, while the partition indexing of dominant Shi regions feels very natural because of the affine symmetric group action, the abaci associated to the partitions do not match up well to the obvious $S_{n}$ action on $\mathbb{Z}^{n}_{n+1}/H$. As we will see, our level $t$ action on $C_{n}^{(mn+1)}$ is an appropriate substitute.

For the level $t$ action, we have:

\begin{lemma} (cf. \cite{F2}, Lemma 3.12) $|C_{n}^{(t)}| = t^{n-1}$ \end{lemma}

\begin{lemma} (cf. \cite{F2}, Proposition 3.13) Each equivalence class in $C_{n}$ under $\equiv_{t}$ contains a unique element of $C_{n}^{(t)}$. \end{lemma}

\begin{lemma} (cf. \cite{F2}, Proposition 3.14) The equivalence relation $\equiv_{t}$ is preserved by the level $t$ action of $\widetilde{S}_{n}$ \end{lemma}

These lemmata together mean in particular that the level $mn \pm 1$ actions are well-defined on the sets $C_{n}^{(mn \pm 1)}$, which have size $(mn \pm 1)^{n-1}$. These are precisely the number of regions and bounded regions in the $m$-Shi arrangement. Our next goal is to understand the orbits of this action. 

\begin{lemma} \label{lem:GAct} The level $mn \pm 1$ action of $\widetilde{S}_{n}$ on $C_{n}^{(t)}$ is generated by the action of $G$. \end{lemma}

Proof. We claim that $s_{0}$ acts on $C_{n}^{(t)}$ in the same way as $s_{\theta}$ does. Since the action of $\widetilde{S}_{n}$ on $\mathbb{Z}$ is completely determined by its action on a set of congruence class representatives mod $n$, we look at the set $Y = \{0,t,2t, \ldots (n-1)t\}$ if $t = mn+1$ and $Y = \{0,-t, -2t, \ldots -(n-1)t\}$ for $t = mn-1$. Then for $1 \leq i \leq n-1$, the action of $s_{i}$ switches the entries of $Y$ congruent to $i-1$ and $i$ mod $n$, and leaves the other elements of $Y$ fixed, so that the action of $s_{\theta}$ on $Y$ switches the elements congruent to 0 and $(n-1)$ mod $n$, and leaves the other elements fixed. But by the same token, for $t=mn+1$, $s_{0} \cdot_{t} 0 = -t \equiv (n-1)t$ mod $nt$, and $s_{0} \cdot_{t} (n-1)t = nt \equiv 0$ mod $nt$, while $s_{0}$ fixes the other elements of $Y$. Similar statements hold for $t = mn-1$, so that, mod $nt$, $s_{0}$ and $s_{\theta}$ have the same level $t$ action on $\mathbb{Z}$, so they act the same way on the equivalence classes under $\sim_{nt}$. $\square$

As a result of this lemma, from this point on, we can focus on the level $t$ action of $G$ on $C_{n}^{(t)}$. To make notation easier, we let $w_{i} = s_{i}$ for $1 \leq i \leq n-1$, and let $w_{0} = s_{\theta}$. Now, following the proof of (\cite{F2}, Lemma 3.10), we describe the stabilizer of the $G$-action on $C_{n}^{(mn+1)}$.

\begin{lemma} \label{lem:Stab} Let $\lambda$ be an $(n,t)$-core partition where $t = mn \pm 1$. Then the $G$-stabilizer of $\lambda$ in $C_{n}^{(t)}$ is generated by the $w_{i}$ for $i$ such that $S(\lambda)_{i} - S(\lambda)_{i-1} =  \mp t$. \end{lemma}

Proof. We prove the case $t = mn+1$ - the other case is similar. If $S(\lambda)_{i} - S(\lambda)_{i-1} = t$, then $w_{i} \cdot_{t} S(\lambda)_{i} = S(\lambda)_{i-1}$ and $w_{i} \cdot_{t} S(\lambda)_{i-1} = S(\lambda)_{i}$, but $w_{i}$ fixes the other elements of $S(\lambda)$, so that $w_{i} \cdot_{t} \lambda = \lambda$.

Now, let $w \in G$ and assume that $w \cdot_{t} \lambda = \lambda$. Since $\lambda$ is a $t$-core, if the abacus for $\lambda$ contains circled beads at $a$ and $a+kt$ for some positive $k$, then the beads $a+t, a+2t, \ldots, a+(k-1)t$ are also circled. The same must be true of the balanced abacus for $\lambda$, so that if $a$ and $a+kt$ are elements of $S(\lambda)$, then so are $a+t, a+2t, \ldots, a + (k-1)t$. Since $t \equiv 1$ mod $n$, we may assume the elements of $S(\lambda)$ lying in a single mod $t$ congruence class are exactly $S(\lambda)_{i} = a$ through $S(\lambda)_{i+k} = a + (k-1)t$.  To fix $S(\lambda)$ setwise, the level $t$ action of $w$ must permute the set $\{S(\lambda)_{i}, S(\lambda)_{i+1}, \ldots, S(\lambda)_{i+k}\}$. But the permutations of this set are generated by $w_{i+1}, w_{i+2}, \ldots, w_{i+k}$. Since the different mod $t$ congruence classes correspond to disjoint cycles in $S_{n}$, the result follows. $\qed$

Of course, the orbit of an $(s,t)$-core $\lambda$ in $C_{n}^{(mn \pm 1)}$ is in bijection with the cosets of its stablizer. This (along with Lemma \ref{floorLemma} below) completes the information we need about the partition side of our main bijection.

\section{Alcoves}

While $n$-sets contain all the information we need about $n$-core partitions, $n$-windows corresponding to alcoves will also be important for us since the ordering of the elements of the window captures the data about the labels of the floors and ceilings of an alcove, as well as the actual hyperplanes that contain those floors and ceilings. For a minimal length coset representative $w$, the permutation $\sigma = f_{n}(w)$ contains the necessary information that links all this data.

To understand the significance of $\sigma$, let $\lambda = w^{-1} \emptyset$. We let $\sigma$ act on $X(\lambda)$ by letting each $w_{i}$ swap the positions of the entries of the tuple that are congruent to $i - 1$ and $i$ mod $n$, which, reduced mod $n$, is the same as the level one action of $s_{i}$. Since $w^{-1} \cdot_{1} (0,1,\ldots, n-1) = X(\lambda)$, the action of $\sigma$ sorts $X(\lambda)$ back in congruence class order $(0,1,2,\ldots)$. As a result, $X(\lambda)_{\sigma(i)} = S(\lambda)_{i-1}.$ (The index shifts by 1 since the ``first'' element of $S(\lambda)$ is $S(\lambda)_{0}$.)

\begin{lemma} Let $w \in \widetilde{S}_{n}$ and let $\alpha_{i}$ be a simple root. Then if $w(\alpha_{i}) = \alpha + k\delta$, $(f_{n}(w))(\alpha_{i}) = \alpha$.
\end{lemma} 

Proof. For $1 \leq i \leq n-1$, $s_{i}$ acts linearly on $Q$ and fixes $\delta$. For $s_{0}$, we note that for any $\alpha'$, $s_{0}(\alpha' + k\delta) = \alpha' + k\delta - \langle \alpha' + k\delta , \alpha_{0} \rangle \alpha_{0} = \alpha' - \langle \alpha' , - \theta \rangle (-\theta) + j\delta$, for some $j$, while $s_{\theta}(\alpha') = \alpha' - \langle \alpha' , \theta \rangle \theta$. Thus the real parts of the actions match. $\qed$

Most importantly, though, $\sigma$ links the labels and the walls of an alcove of $A_{n}$.

\begin{lemma}\label{linkLemma} Let $wA_{0}$ be a dominant alcove, with $\sigma = f_{n}(w)$. Let $H_{\alpha,k}$ be a wall of $wA_{0}$. Then $\sigma s_{i} \sigma^{-1} = s_{\alpha}$ iff $H_{\alpha,k}$ has label $i$. \end{lemma}

Proof. By Lemma \ref{lem:mainFV}, if $H_{\alpha,k}$ is a floor of $wA_{0}$ with label $i$, then $w(\alpha_{i}) = \alpha - k\delta$, so that $\sigma(\alpha_{i}) = \alpha$. If $H_{\alpha,k}$ is a ceiling of $wA_{0}$, then $ws_{i}(\alpha_{i}) = \alpha - k\delta$, so that $\sigma(-\alpha_{i}) = \alpha$. Either way, $(\sigma s_{i} \sigma^{-1}) \cdot \alpha = -\alpha$, so it is a reflection that inverts $\alpha$, and must equal $s_{\alpha}$. $\qed$

So we see that $\sigma$ links $n$-sets with $n$-windows, but also walls of alcoves with their labels. Combining these observations, we have the following lemma about how the partition associated to an $m$-minimal dominant alcove encodes its floors.

\begin{lemma} \label{floorLemma} Let $wA_{0}$ be a dominant alcove that is a subset of $H_{\alpha,k}^{+} \cap H_{\alpha,k+1}^{-}$, where $\alpha = \epsilon_{i} - \epsilon_{j}$ and $i < j$. We let $X$ be the $n$-window of $w^{-1}$ and let $\lambda$ be the partition $w^{-1}\emptyset$ so that $S(\lambda)$ is its $n$-set. Then: 
\begin{enumerate}[i.]
\item $nk < X_{j} - X_{i} < n(k+1)$.
\item $X_{j} - X_{i} = nk+1$ iff $H_{\alpha,k}$ is a wall of $wA_{0}$. 
\item $X_{j} - X_{i} = n(k+1)-1$ iff $H_{\alpha,k+1}$ is a wall of $wA_{0}$.
\item  $|S(\lambda)_{j-1} - S(\lambda)_{i-1}| = X_{\sigma(j)} - X_{\sigma(i)}$.
\end{enumerate}

\end{lemma}

Proof. Part i is intrinsic in \cite{FTV}. In their notation, the alcove $wA_{0}$ has $k_{\alpha} = k$. The numbers $p_{i}$ constructed in section 2.7 of \cite{FTV} are the $n$-window for $w^{-1}$, shifted by a constant. Thus the difference between entries remain the same and so $\lfloor \frac{X_{j} - X_{i}}{n} \rfloor = k$. $X_{j} - X_{i}$ cannot equal $kn$ since $X_{i}$ and $X_{j}$ are not congruent mod $n$.

For the other parts, we note that if $H_{\alpha,k}$ is a wall of $wA_{0}$ with label $r$, then for $\sigma = f_{n}(w)$, we have $\sigma s_{r} \sigma^{-1} = s_{\alpha}$. This means $\sigma$ sends $r$ and $r+1$ to $i$ and $j$ in some order. Thus $X_{j}$ and $X_{i}$ are equal to $S(\lambda)_{r-1}$ and $S(\lambda)_{r}$ in some order. But these differ by $1$ mod $n$ so that $X_{j}-X_{i} = nk+1$ or $n(k+1)-1$. But if $H_{\alpha,k}$ is a wall of $wA_{0}$, then $ws_{r}A_{0}$ is on the other side of $H_{\alpha,k}$, which can only happen if $(s_{r}X)_{j} - (s_{r}X)_{i} < kn$. But since $s_{r}$ only changes the elements of $X$ by a total of 2, this must mean that $X_{j} - X_{i} = kn+1$ and $(s_{r}X)_{j} - (s_{r}X)_{i} = kn-1$. A similar argument proves iii. $\qed$

The conditions ii. and iii. essentially identify floors and ceilings of $wA_{0}$, although not exactly. Hyperplanes of the form $H_{\alpha, 0}$ are the only exception since they are technically never floors or ceilings. And condition iv. completes Lemma \ref{lem:Stab}, showing that it is the labels of the walls of an alcove that index its stabilizer. (Note that it has occurred before that the labels on floors of an alcove index the stabilizer of an object indexing the regions of $\mathcal{S}_{n,m}$. See, for example, Lemma 4.1 in \cite{A2}.)

In what follows, we will want $S_{n}$ to act on $n$-windows of $n$-core partitions, but with the goal of swapping the congruence classes in certain \textit{positions} of the window. To swap, say, positions $i$ and $i+1$ of an $n$-window, we must ``twist'' the action of $w_{i}$ by exactly $\sigma$. More precisely, we have the following:

\begin{lemma} \label{lem:twistAct} Fix $w \in \widetilde{S}_{n}/G$, and $X$, the $n$-window of $w^{-1}$. Let $\sigma = f_{n}(w)$. Then for any $1 \leq i \leq n-1$, $\left((\sigma w_{i} \sigma^{-1}) \cdot_{mn \pm 1} X\right)$ modulo $n$ is exactly $X$ modulo $n$ but with the $i$th and $i+1$st entries swapped. \end{lemma}

Proof. Note that the congruence classes of the $i$th and $i+1$st entries of $X$ are $\sigma(i)$ and $\sigma(i+1)$. But $(\sigma w_{i} \sigma^{-1}) \cdot_{t} \sigma(i)$ is in the congruence class of $\sigma(i+1)$ and vice versa. $\qed$

Notice that, for different choices of $w \in \widetilde{S}_{n}/G$, corresponding to different $(n,mn\pm 1)$-core partitions, the twisting may act differently. So this action is not globally consistent in some sense, but it is well-defined within each orbit of the level $t$ action of $G$ on $C_{n}^{(t)}$.

\section{The Bijections}

We are finally in a position to state our main results.

\begin{theorem} \label{mainThm} Let $w = g_{w}y_{w}$ where $g_{w} \in G$ and $y_{w} \in \widetilde{S}_{n}/G$. Write $\sigma = f(y_{w}^{-1}) \in S_{n}$ and let $X$ be the $n$-window of $y_{w}^{-1}$. The partition associated to $w$ is the partition in $C_{n}^{(mn+1)}$ with $n$-set equivalent to $\sigma g_{w} \sigma^{-1} \cdot_{mn+1} X$ (considered as a set).

This restricts to a bijection between $m$-minimal alcoves $wA_{0}$ and $C_{n}^{(mn+1)}$.

\end{theorem}

\begin{theorem} \label{maxThm} Let $w = g_{w}y_{w}$ where $g_{w} \in G$ and $y_{w} \in \widetilde{S}_{n}/G$. Write $\sigma = f(y_{w}^{-1}) \in S_{n}$ and let $X$ be the $n$-window of $y_{w}^{-1}$. The partition associated to $w$ is the partition in $C_{n}^{(mn-1)}$ with $n$-set equivalent to $\sigma g_{w} \sigma^{-1} \cdot_{mn-1} X$ (considered as a set).

This restricts to a bijection between $m$-maximal alcoves $wA_{0}$ and $C_{n}^{(mn-1)}$.
\end{theorem}

Notice if we assume $g_{w} = e$, then these are precisely the correspondences described in \cite{FV1} and \cite{FV2}. The key result of this paper is that this extends the Fishel-Vazirani results to give a bijection between $C_{n}^{(mn+1)}$ and the set of \textit{all} $m$-minimal alcoves in the $m$-Shi arrangement and then a bijection between $C_{n}^{(mn-1)}$ and the $m$-maximal alcoves. The proofs are virtually the same for both results, merely swapping minimal for maximal, ceilings for floors, and $t=mn+1$ for $t=mn-1$. For clarity, the \textit{exposition} in what follows will discuss the minimal case only, but the \textit{results} will be spelled out for both cases.

A broad perspective on the bijection as a whole might shed some light on the ``twist'' by $\sigma$ in this action. We start with the correspondence between dominant regions and $n$-sets of $(n,mn+1)$-core partitions. The $m$-minimal alcoves can be divided into subsets based on which of them are in the $S_{n}$-orbit of one another. (The $S_{n}$ action restricted to only minimal alcoves is not well-defined, but as we shall see, a relationship like the orbit-stabilizer theorem does exist.) The level $t$ action does the same with the orbits $C_{n}^{(mn+1)}$, a set of the same size as the number of $m$-minimal alcoves. We hope for a bijection between these sets, which we will obtain by way of cosets of certain subgroups on either side of the correspondence. On the partition side, we have the tools in place - the stabilizer (and thus the orbit) under the level $t$ action of $G$ is controlled by the \textit{labels} on the floors of the corresponding alcove. On the alcove side, the group is built from the \textit{walls} of the dominant $m$-minimal alcove. The twisting as described above creates the link between the label of a facet and the root indexing the corresponding wall which will bridge the gap between the two subgroups. Lemma \ref{linkLemma} above merely describes how that link directly relates to the orientation of $wA_{0}$ as related to $A_{0}$. The added relation that comes from substituting $w_{0}$ for $s_{0}$ is a reflection of the fact that $s_{0}A_{0}$ and $w_{0}A_{0}$ have the same orientation.

\begin{figure}[ht]
\label{pic:ex1}
\centerline{\includegraphics[width=4in,height=3.5in]{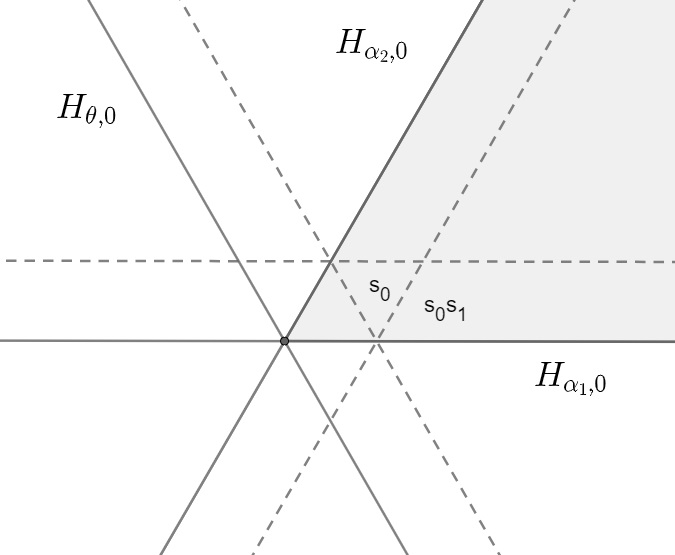}}
\caption{The Shi arrangement $\mathcal{S}_{2,1}$. The shaded area is the dominant chamber.}
\end{figure}

We also give an example illustrating the twist. First we examine the alcove A marked in Figure \ref{pic:ex1} by $w = s_{0}s_{1}$. The 3-window of $w^{-1}$ is $X = (-1,0,4)$, so that $\sigma = s_{2}s_{1}$. Notice that the alcoves $s_{1}A$ and $s_{2}s_{1}A$ are also minimal alcoves, and $\sigma s_{1} \sigma^{-1} = s_{2}s_{1}s_{2}$, while $\sigma s_{2}s_{1} \sigma^{-1} = s_{2}s_{1}$. Then $s_{1}s_{2}s_{1} \cdot_{t} \{0,4,-1\} \equiv_{t} \{3,4,-4\}$, while $s_{2}s_{1} \cdot_{t} \{0,4,-1\} = \{0,7,-4\}$ - these $3$-sets are the only other elements of the orbit of $\{0,4,-1\}$ in $C_{3}^{(4)}$. 

Now, notice that these last calculations could perhaps have been reached more easily - $s_{1}$ could act in a way that mimics the level $t$ action of $\widetilde{S}_{n}$ on $(0,1,2)$ - adding 4 to the first entry, subtracting 4 from the second, and swapping the resulting numbers. Applying this to $\{-1,0,4\}$ (in that order) yields $\{-4,3,4\}$, and then carrying out a similar process on the second and third numbers (mimicing the action of $s_{2}$) gives $\{-4,0,7\}$. However, if we let $A$ be the alcove marked $s_{0}$, it has $n$-window $(-1,1,3)$, and the alcoves $s_{1}A$ and $s_{2}A$ are also minimal. Now, if $s_{1}$ were to shift the first two entries by 4 and then swap them as before, this would yield $(-3,3,3)$, which is no-longer a 3-window! To ``properly'' shift and swap 1 and -1, we must use $\sigma = s_{1}s_{2}s_{1}$ and note that $\sigma s_{1} \sigma^{-1} = s_{2}$. So we apply the level $t$ action of $s_{2}$ to $\{-1,1,3\}$ so that the action adds 4 to 1 and subtracts from -1, properly swapping the mod 3 congruence classes and giving $\{-5,5,3\}$. 

In short, when we write $w = g_{w}y_{w}$, we let $y_{w}^{-1}$ act via the level $1$ action on congruence classes in the $n$-window, and then let $g_{w}$ act via the level $t$ action on positions in the $n$-window. The resulting window indexes the desired partition.

To work toward proving our main result, we first need to show that for every $m$-minimal alcove $wA_{0}$, once we write $w = g_{w}y_{w}$ as above, $y_{w}A_{0}$ is an $m$-minimal alcove.

\begin{lemma} Let $wA_{0}$ be an $m$-minimal alcove of $\mathcal{S}_{n,m}$, and let $\pi \in S_{n}$ be the unique element so $\pi w A_{0}$ is dominant. Then $\pi w A_{0}$ is $m$-minimal.

\end{lemma}

Proof. The action of $S_{n}$ permutes the hyperplanes $H_{\alpha,k}$ for fixed $k$ because the $s_{i}$ for $i>0$ are isometries that fix the origin. Note that if $\pi w A_{0}$ were not $m$-minimal, there would be a hyperplane $H_{\alpha,k}$ with $k >m$ that is a floor of $\pi w A_{0}$, but where that hyperplane is a ceiling of some chamber $\pi w s_{i} A_{0}$. But then $\pi^{-1}(H_{\alpha,k})$ is a floor of $wA_{0}$, a contradiction since $k > m$. $\qed$

This lemma also holds for maximal alcoves - the proof is the same, simply swapping the roles of floors and ceilings.

\begin{lemma} Let $wA_{0}$ be an $m$-maximal alcove of $\mathcal{S}_{n,m}$, and let $\pi \in S_{n}$ be the unique element so $\pi w A_{0}$ is dominant. Then $\pi w A_{0}$ is $m$-maximal.

\end{lemma}

We would like to make an observation about which alcoves in the $S_{n}$-orbit of a dominant $m$-minimal alcove can themselves be $m$-minimal. Specifically, we want to show that these alcoves are also in bijection with the cosets of a certain subgroup of $S_{n}$ built from the dominant alcove. If these subgroups (and thus sets of cosets) correspond to the level $t$ stabilizers of the corresponding partition, then our results are proven. 

We recall that for $I \subseteq \{1,2,\ldots, n-1\}$, the parabolic subgroup $G_{I}$ of $S_{n}$ is the group generated by $\{ (i,i+1) \, | \, i \in I\}$. There is a nice characterization of $G_{I}$'s cosets - specifically, the elements of $G$ that do not have any $(i,i+1)$ for $i \in I$ as a right descent form a set of minimal length left coset representatives for $G_{I}$. (\cite{BB}, section 2.4) We wish to generalize this result to groups generated by \textit{any} set of transpositions in $S_{n}$. First, we note that if $X$ is a set of transpositions in $S_{n}$, then the group $G_{X}$ it generates is a direct product of symmetric groups. Specifically, we define a set partition on $\{1,2,\ldots,n\}$ by choosing the finest set partition so that $j$ and $k$ are in the same set for each $(j,k) \in X$. Then $G_{X}$ is the product of the symmetric groups on each set in this partition.  

\begin{lemma} Let $X$ be a set of transpositions in $S_{n}$. Then the subgroup $G_{X}$ generated by $X$ has a unique set of minimal length left coset reprsentatives in $S_{n}$. \end{lemma}

Proof. Let $w \in G$ and consider the one-line notation for $w$,  $a_{1}a_{2}\ldots a_{n}$ with $a_{i}  = w(i)$ for all $i$. Then if $(j,k)$ is a transposition, the one-line word for $w \cdot (j,k)$ is the same word but with the positions of $j$ and $k$ swapped. Then since $G_{X}$ is a product of full symmetric groups on some sets $X_{i}$, we can rearrange each such set $X_{i}$ by swapping elements of $X_{i}$ within the one-line word, and the resulting word will correspond to another element in the same coset $wG_{X}$. However, there is a unique such word that has each set $X_{i}$ sorted in increasing order, which will srictly minimize the number of inversions. Thus, that one-line word is the unique minimal length coset representative for $wG_{I}. \qed $

Now, we need a characterization of these minimal coset representatives for $G_{X}$. Unfortunately, in general, it is not sufficient to find the elements $w$ so that $\ell(wx) > \ell(w)$ for all $x \in X$. (This is the usual condition for coset representatives of parabolic subgroups. \cite{BB}) Say, for example, that $n=4$ and $X = \{(3,4),(1,4)\} = \{w_{3}, w_{0}\}$. Then the unique minimal set of left coset representatives for $G_{X}$ is $\{e,w_{1}, w_{2}, w_{3}w_{2}\}$. However, the classical definition does not identify these coset representatives in this case. For example, $ w = w_{1}w_{2}w_{1}$ (3214 in one-line notation) has the classical property: right-multiplying by either $(3,4)$ or $(1,4)$ will increase the number of inversions. We must multiply by $(1,3) \notin X$ to get to the minimal coset representative $e$. In some sense, 3214 is a local minimum of the length function within its coset, but not a global minimum. Even though $G_{X}$ is conjugate to a parabolic subgroup of $S_{4}$, the classical coset characterization fails. (It seems that in general, the minimal coset representatives of $G_{X}$ must be described as the $w$ so that $\ell(wy) > \ell(w)$ for all $y \in G_{X}$, so that non-minimal coset representatives come equipped with a weaker property.)

The issue arose here because 1,3, and 4 can be permuted freely, but only using the transpositions $(1,4)$ and $(3,4)$, which are somewhat out of order - this problem would not arise if $X = \{(1,3),(3,4)\}$. Luckily for us, the particular groups $G_{X}$ that arise from the Shi arrangement avoid this obstruction and do follow the classical characterization of minimal coset reprsentatives, which we show now. Consider an $m$-minimal dominant alcove $wA_{0}$. It will have floors of the form $H_{\alpha,k}$ for some positive roots $\alpha$ and $k \leq m$. For each $\alpha$ so that $H_{\alpha,m}$ is a floor of the alcove, we consider the element $s_{\alpha}$ of $G$, considered as the Weyl group of type $A_{n}$. Specifically, if $\alpha = \epsilon_{i} - \epsilon_{j}$ for $i < j$, then $s_{\alpha}$ is the transposition $(i,j)$. We let $X$ be the set of such $s_{\alpha}$ and we will soon show that $G_{X}$ determines which alcoves in the $G$-orbit of $wA_{0}$ are again $m$-minimal. First though, we characterize the minimal coset representatives of these groups.

\begin{lemma} Let $wA_{0}$ be a dominant $m$-minimal alcove of $S_{n,m}$. 
\begin{enumerate} 
 
\item If $H_{\epsilon_{i} - \epsilon_{j},m}$ is a floor of $wA_{0}$ for some $i < j$, then there is no other $k < j$ so that $H_{\epsilon_{k} - \epsilon_{j},m}$ is  a floor of $wA_{0}$.
\item If $H_{\theta,m} $ is a floor of $wA_{0}$, then it is the only floor of the form $H_{\alpha,m}$.
\item If $H_{\epsilon_{i} - \epsilon_{j},m}$ is a ceiling of $wA_{0}$ for some $i < j$, then there is no other $k < j$ so that $H_{\epsilon_{k} - \epsilon_{j},m}$ is a ceiling of $wA_{0}$.
\item If $H_{\theta,m} $ is a ceiling of $wA_{0}$, then it is the only ceiling of the form $H_{\alpha,m}$.

\end{enumerate}
\end{lemma}

Proof. Let $X$ be the $n$-window of $w^{-1}$. First, if $H_{\epsilon_{i} - \epsilon_{j},m}$ is a floor of $wA_{0}$, then the $i$th and $j$th entries of $X$ differ by exactly $mn+1$. But the entries of $X$ are strictly increasing, so if another entry of $X$ differed from the $j$th by $mn+1$, it would come after the $j$th entry.

Then if $H_{\theta,m}$ is a floor of $wA_{0}$, the first and last entries of $X$ differ by exactly $mn+1$. But since $X$ is sorted in ascending order, no other pair of entries can differ by exactly $m$. The last two results again simply swap ceilings for floors and $mn-1$ for $mn+1$. $\qed$

This lemma ensures that the set $X$ cannot be ``out of order'' in the way that caused the problem in our example above if it comes from floors or ceilings of the form $H_{\alpha,m}$ for a particular alcove.

\begin{proposition} Let $wA_{0}$ be a dominant $m$-minimal alcove of $S_{n,m}$. Let $X$ be a set of transpositions which has the property that if $s_{\epsilon_{i}-\epsilon_{j}} \in X$ for $i<j$, then for $k < j$, $k \neq i$, the reflection $s_{\epsilon_{k}-\epsilon_{j}} \notin X$. Define $G^{X} = \{ w \in S_{n} \, | \, \ell(wx) > \ell(w) \textrm{ for all } x \in X\}$. Then $G^{X}$ is the set of minimal coset representatives for $G_{X}$.

\end{proposition}

Proof. First, it is clear that any minimal coset representative for $G_{X}$ has the property that $\ell(wx) > \ell(w)$ for all $x \in X$. We must merely show that, for any $w \in S_{n}$, if there is some $y \in G_{X}$ with $\ell(wy) < \ell(w)$, then in fact there is an $x \in X$ with $\ell(wx) < \ell(w)$ so that $w \notin G^{X}$.

As before, we consider $a_{1}a_{2}\ldots a_{n}$, the word for $w$ in one-line notation. If $\ell(wy) < \ell(w)$ for some $y \in G_{X}$, then choose some inversion $a_{i} > a_{j}$ with $i < j$ so that $a_{i}$ occurs after $a_{j}$ in the word for $wy$. This means that the set partition of $\{1,2, \ldots, n\}$ associated to $X$ contains some set $B$ that contains both $a_{i}$ and $a_{j}$. Write $B = (b_{1},b_{2}, \ldots b_{k})$ in increasing order. Then by the previous lemma, $X$ must contain the transpositions $(b_{1},b_{2})$, $(b_{2}, b_{3}), \ldots, (b_{k-1},b_{k})$. If $a_{i}$ and $a_{j}$ appear consecutively in $B$, then $(a_{i},a_{j}) \in X$, and $\ell(w \cdot (a_{i},a_{j})) < \ell(w)$, as desired. Otherwise, we proceed inductively. 

If they are not consecutive in $B$, choose some $a_{r}$ with $a_{j} < a_{r} < a_{i}$. If $r < j$, then $a_{r}$ and $a_{j}$ are an inversion in $w$. If $r > j$, then $r > i$, so that $a_{r}$ and $a_{i}$ are an inversion in $w$. Either way, we have found another pair of elements of $B$ which form an inversion in $w$ and whose entries in $B$ are closer together than $a_{i}$ and $a_{j}$. Repeating this, we eventually find consecutive elements of $B$ that are inverted in $w$. $\qed$

This characterization of the groups $G_{X}$ finally lets us show that the $m$-minimal alcoves in the $S_{n}$ orbit of an $m$-minimal alcove are in bijection with $G^{X}$.

\begin{theorem} \label{orbitThm} Let $w \in \widetilde{S}_{n}$ such that $wA_{0}$ is dominant and $m$-minimal. Let $X$ be the set of transpositions $s_{\alpha}$ so that $H_{\alpha,m}$ is a floor of the alcove $wA_{0}$. Let $\rho \in S_{n}$. Then $\rho w A_{0}$ is $m$-minimal if and only if $\rho \in G^{X}$. \end{theorem}

Proof. First, note that the action of $G$ permutes the hyperplanes $H_{\alpha,k}$ for fixed $k$. Thus the only way that $\rho wA_{0}$ can fail to be $m$-minimal is if $\rho$ sends a hyperplane $H_{\alpha,m}$ that is a floor of $w A_{0}$ to a hyperplane $H_{\alpha',m}$ which is a floor of $\rho w A_{0}$, where $\alpha' < 0$. By Lemma \ref{lem:mainFV}, if $H_{\alpha,m}$ is a floor of $wA_{0}$ with label $j$, then $w(\alpha_{j}) = \alpha-m\delta$. But then $\rho w(\alpha_{j}) = \rho(\alpha) - m\delta$. So we let $\alpha' = \rho(\alpha)$ and note that if $\alpha' < 0$, then $\alpha$ is an inversion of $\rho$, so that $\ell(\rho s_{\alpha}) < \ell(\rho)$. This cannot occur if $\rho \in G^{X}$. 

Conversely, though, if $\rho \notin G^{X}$, there must be some $s_{\alpha} \in X$ so that $\ell(\rho s_{\alpha}) < \ell(\rho)$. Assume that the facet of $wA_{0}$ contained in $H_{\alpha,m}$ has label $j$. Then $w(\alpha_{j}) = \alpha - m\delta$, and by assumption $\rho w (\alpha_{j}) = \alpha' - m\delta$ where $\alpha' < 0$. This implies that $H_{\alpha',m}$ is a floor of $\rho w A_{0}$, which is thus not $m$-minimal. $\qed$

This proof needs very few adjustments to establish the following.

\begin{theorem} \label{orbitThm2} Let $w \in \widetilde{S}_{n}$ such that $wA_{0}$ is dominant and $m$-maximal. Let $X$ be the set of transpositions $s_{\alpha}$ so that $H_{\alpha,m}$ is a ceiling of the alcove $wA_{0}$. Let $\rho \in S_{n}$. Then $\rho w A_{0}$ is $m$-minimal if and only if $\rho \in G^{X}$. \end{theorem}

Proof. Now, the only way that $\rho wA_{0}$ can fail to be $m$-maximal is if $\rho$ sends a hyperplane $H_{\alpha,m}$ that is a ceilng of $w A_{0}$ to a hyperplane $H_{\alpha',m}$ which is a ceiling of $\rho w A_{0}$, where $\alpha' < 0$. By Lemma \ref{lem:mainFV}, if $H_{\alpha,m}$ is a ceiling of $wA_{0}$ with label $j$, then $w(\alpha_{j}) = m\delta - \alpha$. But then $\rho w(\alpha_{j}) = m\delta - \rho(\alpha)$. So we let $\alpha' = \rho(\alpha)$ and note that if $\alpha' < 0$, then $\alpha$ is an inversion of $\rho$, so that $\ell(\rho s_{\alpha}) < \ell(\rho)$. This cannot occur if $\rho \in G^{X}$. 

Conversely, though, if $\rho \notin G^{X}$, there must be some $s_{\alpha} \in X$ so that $\ell(\rho s_{\alpha}) < \ell(\rho)$. Assume that the facet of $wA_{0}$ contained in $H_{\alpha,m}$ has label $j$. Then $w(\alpha_{j}) = m\delta -\alpha$, and by assumption $\rho w (\alpha_{j}) = m\delta - \alpha'$ where $\alpha' < 0$. This implies that $H_{\alpha',m}$ is a ceiling of $\rho w A_{0}$, which is thus not $m$-maximal. $\qed$

We now understand the alcove side of the bijection, which essentially completes the proof of our two main results.

Proof of Theorem \ref{mainThm} and \ref{maxThm}. The $m$-minimal alcoves of $\mathcal{S}_{n,m}$ are in bijection with pairs consisting of a dominant $m$-minimal alcove $wA_{0}$ and a minimal coset representative of $G_{X}$ as defined in Theorem $\ref{orbitThm}$. For any such pair, by Lemmas \ref{linkLemma} and \ref{lem:Stab}, $G_{X}$ is conjugate to $H$, the $G$-stabilizer of the partition $w^{-1} \emptyset$ in $C_{n}^{(mn+1)}$ under the level $t$ action. Thus the cosets of $G_{X}$ and $H$ are in bijection.

Similarly, the $m$-maximal alcoves of $\mathcal{S}_{n,m}$ are in bijection with pairs consisting of a dominant $m$-maximal alcove $wA_{0}$ and a minimal coset representative of $G_{X}$ as defined in Theorem $\ref{orbitThm2}$. For any such pair, by Lemmas \ref{linkLemma} and \ref{lem:Stab}, $G_{X}$ is conjugate to $H$, the $G$-stabilizer of the partition $w^{-1} \emptyset$ in $C_{n}^{(mn-1)}$ under the level $t$ action. Thus the cosets of $G_{X}$ and $H$ are in bijection. $\qed$

Before we close, we want to show one relationship between our results and other combinatorial indexings of Shi regions. Specifically, \cite{AL} assigns an $m$-parking function to each region of $\mathcal{S}_{n,m}$. We show that our bijection reflects the structure of these parking functions - specifically, their stabilizers. (Note this argument is essentially in \cite{A2}, but we include it here for completeness.)

\begin{proposition} \label{prop:ALProp} Let $wA_{0}$ be a dominant $m$-minimal alcove, with $X$ the transpositions associated to floors of $wA_{0}$ as defined as in Theorem \ref{orbitThm}. Let $f$ be the $m$-parking function associated to $wA_{0}$ in \cite{AL}. Then if $x \in X$, $x \cdot f = f$ when acting via permutation.
     
\end{proposition}

Proof. The Athanasiadis-Linusson bijection is as follows. Make a list consisting of terms of the form $x_{i} + k$ for $0 \leq k \leq m-1$. Arrange the list in the unique permutation $\{y_{i}\}$ so that for points in $wA_{0}$, $y_{1} > y_{2} > \ldots > y_{mn}$ holds. We then connect each $x_{i} + k$ for $k>0$ to $x_{i} + k-1$ with an arc. In addition, we connect $x_{i}$ to $x_{j}+m-1$ with an arc whenever $wA_{0} \subseteq H_{\epsilon_{i} - \epsilon_{j},m}^{+}$. Finally we remove each arc containing another and define $f(i)$ to be the position of the leftmost term in the chain containing $x_{i}$.

So, assume $s_{\epsilon_{i} - \epsilon_{j}} \in X$ for $i < j$, so that the construction for the parking functions connects $x_{i}$ to $x_{j} + m-1$. Thus, as long as this arc is not deleted in the construction, $f(i) = f(j)$. However, this arc would only be deleted if another arc were contained in it. So assume $x_{i} > x_{r} + k  > x_{s} + \ell > x_{j} + m-1$ is a sublist in the construction with $x_{r} + k$ connected to $x_{s} + \ell$. 

For the first case, this configuration could occur if $r = s$ and $\ell = k-1$. Note that $x_{i} > x_{r} > x_{s}$, so $i > r > s$ by dominance. But then $x_{i} - x_{r} > k$, and $x_{r} - x_{j} > m-k$. If $H_{\epsilon_{i}-\epsilon_{j},m}$ were a floor of the region, then there would be an alcove of $\mathcal{A_{n}}$ sharing that facet with $wA_{0}$. However, $x_{i} - x_{r} > k$, and $x_{r} - x_{j} > m-k$ would still hold on that region, contradicting the fact that $x_{i} - x_{j} < m$.

For the second case, it could occur that $k=0$ while $\ell = m-1$. However, this would mean that $x_{i}  > x_{r} > x_{s} > x_{j}$, so that $i > r > s > j$ by dominance. However, we then have $x_{i} - x_{r} > 0$, $x_{r} - x_{s} > m$, and $x_{s} - x_{j} > 0$. Again, it is impossible for $H_{\epsilon_{i} - \epsilon_{j},m}$ to be a floor of this region by the same argument.

So, $X$ is a subset of the stabilizer of $f$. However, the index of the stabilizer of $f$ is the same as the index of $G_{X}$ by Theorem \ref{mainThm}, so $G_{X}$ has the same number of elements as the stabilizer of $f$ and they must be equal. $\qed$

\section{Conclusion}

Our two main results provide a unified extension of the results of \cite{FV1} and \cite{FV2}. It is worth noting that these bijections reflect other combinatorial indexings of Shi regions. In essence, Proposition \ref{prop:ALProp} shows that  the level $t$ action translates the $S_{n}$-action inherent in the Athanasiadis-Linusson bijection (along with the corresponding stabilizer and orbits) and translates them into the partition context of Fishel-Vazirani. The labelled Dyck paths of Garsia and Haiman (\cite{GH}) are similar - each valid Dyck path corresponds to a dominant $m$-minimal alcove, and the labels are permutations which represent the different Weyl chambers with a minimal alcove in the $S_{n}$ orbit of the dominant one. Whether $C_{n}^{(mn+1)}$ is the exact right indexing set to use is somewhat questionable - the $G$-orbits of the $(n,mn+1)$ cores (ignoring the equivalence $\equiv_{mn+1}$) would work just as well to index the minimal alcoves, but a description of that set of partitions is somewhat elusive. A better understanding of how the level $t$ action acts directly on abaci might be helpful in understanding it. Another area of interest might be understanding how to turn a partition in $C_{n}^{(mn+1)}$ directly into a parking function. For $(n,n+1)$-core partitions, the process is relatively straightforward. We can take the positive abacus corresponding to $\lambda$ and read it as a Dyck path as in \cite{A} - the area above that path forms a partition $\lambda'$ that fits in the ``staircase'' partition $(n-1,n-2,\ldots,1)$. Then the removable boxes in $\lambda'$ encode roots as described in \cite{FTV}. The Athanasiadis-Linusson diagram corresponding to the same region as $\lambda$ is $12\ldots n$, with an arc connecting $i$ and $j$ exactly if the box corresponding to the root $\epsilon_{i}-\epsilon_{j}$ is removable in $\lambda'$. A similar process to go directly from $C_{n}^{(mn+1)}$ to an $m$-parking function would be greatly desirable.

It is also worth noting that near the completion of this paper, the author became aware of \cite{DHFM}, which takes a vastly different approach to characterizing the $m$-minimal alcoves than Theorem \ref{orbitThm}. A better understanding of how these characterizations interact would be greatly desirable - in particular determining the extent to which any of our results can be generalized out of type A seems like an important next step.

\end{document}